\documentclass[12pt,authoryear,preprint]{elsarticle}
\usepackage{pgf}
\usepackage{amsmath}   
\usepackage[ansinew]{inputenc}                                                                                                                                                                                                                
\usepackage{amsmath}
\usepackage{amsmath,amssymb}
\usepackage{multirow}
\usepackage[]{algorithm2e}
\usepackage{rotating}
\usepackage{enumerate}
\usepackage{array}
\usepackage[labelformat=simple]{subfig}
\usepackage{float}
\usepackage{color}
\usepackage{colortbl}
\usepackage{wasysym}
\usepackage{mathrsfs}
\usepackage{enumerate}
\usepackage{pxfonts}
 \usepackage{amsfonts}
\usepackage{dsfont}
\usepackage{array}
\usepackage{lscape}
\usepackage[pdftex,                %
    bookmarks         = true,%     % Signets
    bookmarksnumbered = true,%     % Signets num\'erot\'es
    pdfpagemode       = None,%     % Signets/vignettes ferm\'e Ã  l'ouverture
    pdfstartview      = FitH,%     % La page prend toute la largeur
    pdfpagelayout     = SinglePage,% Vue par page
    colorlinks        = true,%     % Liens en couleur
    urlcolor          = magenta,%  % Couleur des liens externes
    pdfborder         = {0 0 0}%   % Style de bordure : ici, pas de bordure
    ]{hyperref}% 

\usepackage{tabularx,ragged2e,booktabs,caption}
\newcolumntype{Y}{>{\centering\arraybackslash}X}
\newtheorem{The}{Theorem}[section]

\newtheorem{Def}[The]{Definition}

\DeclareMathOperator{\Cov}{Cov}

\usepackage{color}
\definecolor{darkred}{rgb}{.7, 0, 0}

\numberwithin{equation}{section}%num\'eros par section
\numberwithin{table}{section}
\numberwithin{figure}{section}
\usepackage{geometry}

\journal{publication.}

\begin{document}

\begin{frontmatter}
\date{February 5, 2014}
\title{Effects of associated kernels in nonparametric multiple regressions}% \thanks[label1]
\author[rvt]{Sobom M. Som\'e\corref{cor1}}
\ead{sobom.some@univ-fcomte.fr}
%\author[rvt1]{Tristan   Senga Ki\'ess\'e}
% \ead{tristan.sengakiesse@univ-nantes.fr}
\author[rvt]{C\'elestin C. Kokonendji}
\ead{celestin.kokonendji@univ-fcomte.fr}

\cortext[cor1]{\textit{Corresponding author:} Universit\'e de Franche-Comt\'e, UFR Sciences et Techniques, Laboratoire de Math\'ematiques de Besan\c{c}on -- UMR 6623 CNRS-UFC, 16 route de Gray, 25030 Besan\c{c}on cedex, France; Tel. +33 381 666 398; Fax +33 381 666 623.}

\address[rvt]{University of Franche-Comt\'e, Besan\c{c}on, France}
%\address[rvt1]{University of Nantes Angers Le Mans, France}
\begin{abstract}
Associated kernels have been introduced to improve the classical continuous kernels for smoothing any functional on several kinds of supports such as bounded continuous and discrete sets. This work deals with the effects of combined associated kernels on nonparametric multiple regression functions. Using the Nadaraya-Watson estimator with optimal bandwidth matrices selected by cross-validation procedure, different behaviours of multiple regression estimations are pointed out according the type of multivariate associated kernels with correlation or not. Through simulation studies, there are no effect of correlation structures for the continuous regression functions and also for the associated continuous kernels; however, there exist really effects of the choice of multivariate associated kernels following the support of the multiple regression functions bounded continuous or discrete. Applications are made on two real datasets. %A large $R^2$ is not synomym of adjustment quality essentially for not suitable associated kernels.
\end{abstract}

\begin{keyword}
Bandwidth matrix \sep  continuous associated kernel \sep correlation structure \sep cross-validation \sep discrete associated kernel \sep Nadaraya-Watson estimator.\newline
 {\it Mathematics Subject Classification 2010}: 62G05(08); 62H12

%\MSC

{\bf Short Running Title}: Associated kernels in multiple regressions
\end{keyword}
%{\it MSC 2010:} Primary 62H05, 60E07; Secondary 62F10, 62H99.
\end{frontmatter}

\section{Introduction}

Considering the relation between a response variable $Y$ and a $d$-vector $(d\geq1)$ of explanatory variables $\mathbf{x}$ given by
\begin{equation}\label{Regfonction}
Y = m\left(\mathbf{x}\right) + \epsilon,
\end{equation}where $m$ is the unknown regression function from $\mathbb{T}_{d} \subseteq \mathbb{R}^{d}$ to $\mathbb{R}$ and   $\epsilon$ the disturbance term with null mean and finite variance. Let $(\bold{X}_{1}, Y_{1}), \ldots, (\bold{X}_{n}, Y_{n})$ be a sequence of independent and identically distributed (iid) random vectors on $\mathbb{T}_{d} \times \mathbb{R} (\subseteq \mathbb{R}^{d+1})$ with  $m(\bold{x}) =  \mathbb{E}\left(Y|\bold{X} = \bold{x}\right)$ of (\ref{Regfonction}). The  \cite{Nadaraya69} and \cite{Watson64} estimator $\widehat{m}_{n}$ of $m$, using continuous classical (symmetric) kernels  is
\begin{equation}\label{RegClassic}
  \widehat{m}_{n}(\mathbf{x};K,\mathbf{H}) =  \displaystyle \sum_{i = 1}^{n} \dfrac{Y_iK\left\{\mathbf{H}^{-1}(\mathbf{x} - \mathbf{X}_i)\right\}}{ \sum_{i = 1}^{n}K\left\{\mathbf{H}^{-1}(\mathbf{x} - \mathbf{X}_i)\right\}}=  \widehat{m}_{n}(\mathbf{x};\mathbf{H}) , ~~\forall \mathbf{x} \in   \mathbb{T}_{d} := \mathbb{R}^{d},
\end{equation} where $\mathbf{H}$ is the symmetric and positive definite bandwidth matrix of dimension $d\times d$ and the function $K(\cdot)$ is the multivariate kernel assumed to be spherically symmetric probability density function. Since the choice of the kernel $K$ is not important in classical case, we use the common notation  $\widehat{m}_{n}(\mathbf{x};\mathbf{H})$ for classical kernel regression. The multivariate classical kernel (e.g. Gaussian) suits only for regression functions on unbounded supports (i.e. $\mathbb{R}^d$); see also \cite{S92}.  \cite{Racine&Li04} proposed product of  kernels composed by univariate Gaussian kernels for continuous variables and \cite{AJG76} kernels for categorical variables; see also \cite{Hayfield&Racine} for some implementations and uses of these multiple kernels. Notice that the use of symmetric kernels gives weights outside  variables with unbounded supports. In the univariate continuous case, \cite{{Chen99},{Chen00b}} is one of the pioneers who has proposed asymmetric kernels (i.e. beta and gamma) which supports coincide with those of the functions to be estimated.  \cite{Zhang2010} and \cite{ZhangandKarunamuni2010} studied the performance of these beta and gamma kernel estimators at the boundaries in comparison with those of the classical kernels. Recently, \cite{L13} investigated several families of these univariate continuous kernels that he called univariate associated kernels; see also \cite{KSKZ07}, \cite{KSK11}, \cite{ZAK12} and \cite{WKK14} for univariate discrete situations. A continuous multivariate version of these associated kernels have been studied by \cite{Kokonendji&Some14} for density estimation.

The main goal of this work is to consider  multivariate associated kernels and then to investigate the importance of their choice   in multiple regression. These associated kernels are appropriated for both continuous and count explanatory variables. In fact, in order to estimate the regression function $m$ in (\ref{Regfonction}), we propose multiple (or product of) associated  kernels composed by univariate discrete associated kernels (e.g. binomial, discrete triangular) and continuous ones (e.g. beta, Epanechnikov). We will  also use a bivariate beta kernel with correlation structure. Another motivation of this work is to investigate the effect of correlation structure for  explanatory variables in  continuous regression estimation. These associated kernels suit for this situation of mixing axes as  they fully respect the support of each explanatory variable. In other words, we will measure the effect of type of associated kernels, denoted $\boldsymbol{\kappa}$, in multiple regression by simulations and applications.

The rest of the paper is organized as follows.  Section~\ref{sec:Multiple regression} gives a general definition of multivariate associated kernels which includes the continuous classical symmetric and the multiple composed by univariate discrete and continuous. For each definition, the corresponding kernel regression  appropriated for both continuous and discrete explanatory variables are given. In Section~\ref{sec:Simulation studies and real data analisys}, we explore  the importance of the choice of appropriated associated kernels according to the support of the variables through simulations studies  and real data analysis.  Finally, summary and final remarks  are drawn in Section~\ref{sec:Summary and final remarks}.

%Furthermore, some comparisons are pointed out within continuous kernels (multiple and general) and also mixed (multiple) associated kernels.
\section{Multiple regression by associated kernels}
\label{sec:Multiple regression}

\subsection{Definition}
\label{sec:Multivariate associated regressions}
In order to include both discrete and continuous regressors, we assume $\mathbb{T}_{d}$ is any subset of $\mathbb{R}^d$. More precisely, for $j = 1, \ldots, n$, let us consider on 
$\mathbb{T}_{d} = \displaystyle \otimes_{j=1}^{d}\mathbb{T}^{[j]}_1$ the measure $\boldsymbol{\nu} = \nu_{1} \otimes \ldots  \otimes \nu_{d}$ where $\nu_{j}$ is a Lesbesgue or count measure on the corresponding univariate support $\mathbb{T}^{[j]}_1$. Under these assumptions, the associated kernel  $K_{\bold{x}, \mathbf{H}}(\cdot)$ which replaces the classical kernel $K(\cdot)$ of (\ref{RegClassic})  is a probability density function (pdf) in relation to a measure $\boldsymbol{\nu}$. This  kernel $K_{\bold{x}, \mathbf{H}}(\cdot)$  can be defined as follows.

\begin{Def}\label{Defkern}
Let $\mathbb{T}_{d}\left(\subseteq \mathbb{R}^{d}\right)$ be the support of the regressors,  $\mathbf{x}  \in \mathbb{T}_{d}$ a target vector and   $\mathbf{H}$ a bandwidth matrix. A parametrized pdf $K_{\mathbf{x}, \mathbf{H}}(\cdot)$ of support $\mathbb{S}_{\mathbf{x}, \mathbf{H}}  \left(\subseteq \mathbb{R}^{d}\right)$
 is called ``multivariate (or general) associated kernel'' if the following conditions are satisfied:
\begin{align}
  &\mathbf{x} \in \mathbb{S}_{\mathbf{x}, \mathbf{H}},  \label{NoyAss1}\\
  &\mathbb{E}\left(\mathcal{Z}_{\mathbf{x}, \mathbf{H}}\right) = \mathbf{x} +  \mathbf{a}(\mathbf{x}, \mathbf{H}), \label{NoyAss2} \\
  & \rm{Cov}\left(\mathcal{Z}_{\mathbf{x}, \mathbf{H}}\right) =   \mathbf{B}(\mathbf{x}, \mathbf{H}), \label{NoyAss3}
\end{align}
where $\mathcal{Z}_{\mathbf{x}, \mathbf{H}}$ denotes the random vector with pdf  $K_{\mathbf{x},  \mathbf{H}}$ and both $\mathbf{a}(\mathbf{x},  \mathbf{H}) = \left(a_{1}(\mathbf{x},  \mathbf{H}), \ldots, a_{d}(\mathbf{x},  \mathbf{H})\right)^{\top}$ and  $\mathbf{B}(\mathbf{x},  \mathbf{H}) = \left(b_{ij}(\mathbf{x},  \mathbf{H})\right)_{i,j = 1, \ldots, d}$ tend, respectively, to the null vector $\mathbf{0}$ and the null matrix $\mathbf{0}_d$ as $\mathbf{H}$ goes to $\mathbf{0}_d$.
\end{Def}
From this definition and in comparison with (\ref{RegClassic}),  the Nadaraya-Watson estimator using associated kernels is 

\begin{equation}\label{RegFull}
\widetilde{m}_{n}(\bold{x};K_{\bold{x}, \mathbf{H}}) = \displaystyle \sum_{i = 1}^{n} \dfrac{Y_iK_{\bold{x}, \mathbf{H}}\left(\bold{X}_i\right)}{ \sum_{i = 1}^{n}K_{\bold{x}, \mathbf{H}}\left(\bold{X}_i\right)}=\widetilde{m}_{n}(\bold{x};\boldsymbol{\kappa},\mathbf{H}), ~~~ \forall  \bold{x} \in  \mathbb{T}_{d} \subseteq \mathbb{R}^{d},
\end{equation} where $\mathbf{H} \equiv \mathbf{H}_{n} $ is the bandwidth matrix such that $ \mathbf{H}_{n} \rightarrow \mathbf{0}$ as $n \rightarrow  \infty$, and $\boldsymbol{\kappa}$ represents the type of associated kernel $K_{\bold{x}, \mathbf{H}}$, parametrized by $\bold{x}$ and $\mathbf{H}$.  Without loss of generality and to point out the effect of $\boldsymbol{\kappa}$, we  will in hereafter use $\widetilde{m}_{n}(\bold{x};\boldsymbol{\kappa},\mathbf{H}) = \widetilde{m}_{n}(\bold{x};\boldsymbol{\kappa})$ since the bandwidth matrix is here investigated only by cross validation.

The following two examples provide the well-known and also interesting particular cases of multivariate associated kernel estimators. The first can be seen as an interpretation of  classical associated kernels through continuous symmetric kernels. The second  deals on non-classical associated kernels without correlation structure.

 Given a target vector $\mathbf{x}\in\mathbb{R}^{d} =:\mathbb{T}_d$ and a bandwidth matrix $\mathbf{H}$, it follows that the classical kernel in (\ref{RegClassic}) with null mean vector and  covariance matrix $\mathbf{\Sigma}$ induces the so-called (multivariate) classical  associated kernel: 
\begin{eqnarray}\label{classcical}
(i)~~K_{\mathbf{x}, \mathbf{H}}(\cdot) = \dfrac{1}{\det\mathbf{H}} K\left\{\mathbf{H}^{-1}(\mathbf{x} - \cdot)\right\}
\end{eqnarray} on $\mathbb{S}_{\mathbf{x}, \mathbf{H}} = \mathbf{x} - \mathbf{H}\mathbb{S}_{d}$ with $\mathbb{E}\left(\mathcal{Z}_{\mathbf{x}, \bf{H}}\right) = \mathbf{x}$ (i.e. $\mathbf{a}(\mathbf{x}, \mathbf{H}) = \mathbf{0}$) and $\rm{Cov}\left(\mathcal{Z}_{\mathbf{x}, \mathbf{H}}\right) =  \mathbf{H}\mathbf{\Sigma}\mathbf{H}$; 
\begin{eqnarray*}
(ii)~~K_{\mathbf{x}, \mathbf{H}}(\cdot) = \dfrac{1}{(\det\mathbf{H})^{1/2}} K\left\{\mathbf{H}^{-1/2}(\mathbf{x} - \cdot)\right\}
\end{eqnarray*} on $\mathbb{S}_{\mathbf{x}, \mathbf{H}} = \mathbf{x} - \mathbf{H}^{1/2}\mathbb{S}_{d}$ with $\mathbb{E}\left(\mathcal{Z}_{\mathbf{x}, \bf{H}}\right) = \mathbf{x}$ (i.e. $\mathbf{a}(\mathbf{x}, \mathbf{H}) = \mathbf{0}$) and $\rm{Cov}\left(\mathcal{Z}_{\mathbf{x}, \mathbf{H}}\right) =  \mathbf{H}^{1/2}\mathbf{\Sigma}\mathbf{H}^{1/2}$.

A second particular case of  Definition~\ref{Defkern}, appropriate for both continuous and count explanatory variables without correlation structure is presented as follows.

Let $\mathbf{x} = (x_{1}, \ldots, x_{d})^{\top} \in \times_{j=1}^{d}\mathbb{T}^{[j]}_{1} =:\mathbb{T}_d$  and $\mathbf{H} = \mathbf{Diag} (h_{11}, \ldots, h_{dd})$ with $h_{jj} > 0$. Let  $K^{[j]}_{x_{j}, h_{jj}}$ be a (discrete or continuous) univariate associated kernel (see Definition~\ref{Defkern} for $d = 1$) with its corresponding random variable $\mathcal{Z}^{[j]}_{x_{j}, h_{jj}}$  on $\mathbb{S}_{x_{j}, h_{jj}} (\subseteq \mathbb{R})$ for all $j = 1, \ldots, d$. Then,  the multiple associated kernel is also a multivariate associated kernel:
%\begin{Pro} \label{Proprod}
%Let $\times_{j=1}^{d}\mathbb{T}^{[j]}_{1} = \mathbb{T}_{d}$ be the support of the density $f$ to be estimated with $\mathbb{T}^{[j]}_{1} (\subseteq \mathbb{R})$ the supports of univariate margins of $f$. 
\begin{eqnarray}\label{prodkern1}
  K_{\mathbf{x}, \mathbf{H}}(\cdot) = \prod_{j=1}^{d}K^{[j]}_{x_{j}, h_{jj}}(\cdot)
\end{eqnarray}on $\mathbb{S}_{\mathbf{x}, \mathbf{H}} =  \displaystyle \times_{j=1}^{d}\mathbb{S}_{x_{j}, h_{jj}}$ with $\mathbb{E}\left(\mathcal{Z}_{\mathbf{x}, \bf{H}}\right) = \left(x_{1} + a_{1}(x_{1}, h_{11}), \ldots, x_{d} + a_{d}(x_{d}, h_{dd})\right)^{\top}$ and $\Cov\left(\mathcal{Z}_{\mathbf{x}, \mathbf{H}}\right)$ = $ \mathbf{Diag}\left(b_{jj}(x_{j}, h_{jj})\right)_{j = 1, \ldots, d}$. In other words, the random variables  $\mathcal{Z}^{[j]}_{x_{j}, h_{jj}}$ are independent components of the random vector $\mathcal{Z}_{\mathbf{x}, \bf{H}}$.

Here, in addition to the Nadaraya-Watson estimator using general associated kernels  given in (\ref{RegFull}), we proposed a slight one. In fact, for multivariate supports composed of continuous and discrete univariate support, we lack appropriate general associated kernels. Therefore, the estimator (\ref{RegFull}) becomes with multiple associated kernels (\ref{prodkern1}):
\begin{equation}
\widetilde{m}_{n}(\bold{x};\boldsymbol{\kappa})= \displaystyle \sum_{i = 1}^{n} \dfrac{Y_i\prod_{j=1}^{d}K^{[j]}_{x_{j}, h_{jj}}(X_{ij})}{ \sum_{i = 1}^{n}\prod_{j=1}^{d}K^{[j]}_{x_{j}, h_{jj}}(X_{ij})},~~~\forall \bold{x} = \left(x_{1}, \ldots, x_{d}\right)^{\top}  \in \mathbb{T}_{d} :=  \displaystyle \times_{j=1}^{d} \mathbb{T}^{[j]}_{1}. \label{Regprod}
\end{equation} In theory and in practice, one often uses (\ref{Regprod}) from multiple associated kernels (\ref{prodkern1}) which are more manageable than (\ref{RegFull}); see, e.g., \cite{S92} and also \cite{BR10}  for density estimation.

%The third particular case of multivariate associated kernels is constructed from a pdf composed of product of univariate pdf and a correlation structure using the correlation structure introduced by .   This kind of associated kernels enables to reach some points of the multidimensional smoothing.  \cite{Kokonendji&Some14} illustrated the effect of this  technique on  bivariate beta kernel and for density estimation.  The Nadaraya-Watson estimator is then using such  associated kernel with correlation structure in (\ref{RegFull}). 
%Like \cite{Bertin&Klutnitchkoff14}, the  miminax properties of this bivariate beta kernel are also possible and more generally for associated kernels.
\subsection{Associated kernels for illustration}
\label{ssec:Associated kernels for illustration}
In order to point out  the importance  of the type of kernel  $\boldsymbol{\kappa}$ in a regression study, we motivate below some kernels that will be used in simulations. These concern seven basic associated kernels for which three of them are univariate discrete,  three others are  univariate continuous and the last one is  a bivariate beta with correlation structure.
\begin{itemize}
\item  The binomial kernel ($\texttt{Bin}$) is defined on the support  $\mathbb{S}_x = \{0, 1,  \ldots, x+1\}$  with $x \in \mathbb{T}_1 :=\mathbb{N}= \{0,1,\ldots\}$ and then $h \in (0, 1]$: 
\begin{equation*}\label{c}
B_{x,h}(u)=\frac{(x+1)!}{u!(x+1-u)!}\left(\frac{x+h}{x+1}\right)^{u}\left(\frac{1-h}{x+1}\right)^{x+1-u} \mathds{1}_{\mathbb{S}_x}(u),
\end{equation*} where $\mathds{1}_{A}$ denote the indicator function of any given event $A$. Note that $B_{x,h}$ is the probability mass function (pmf) of the binomial distribution $\mathcal{B}(x+1; (x+h)/(x+1))$ with its number of trials  $x+1$ and its success probability in each trial $(x+h)/(x+1)$.
It is appropriated for count data with small or moderate sample sizes and, also, it does not  satisfy (\ref{NoyAss3}); see~ \cite{KSK11} and also \cite{ZAK12} for a bandwidth selection by Bayesian method.

\item For fixed arm $a \in \mathbb{N}$, the discrete  triangular kernel ($\texttt{DTr}a$) is defined on  $\mathbb{S}_{x,a} = \left\{x,  x \pm 1, \ldots, x \pm a \right\}$   with $x \in \mathbb{T}_1 =\mathbb{N}$: 
   \begin{equation*}\label{c}
 DT_{x, h;a}(u)=\frac{(a+1)^{h} - |u - a|^{h}}{P(a, h)}\mathds{1}_{\mathbb{S}_{x\setminus \{a\}} }(u), 
\end{equation*}where $h >0$ and  $P(a, h)= (2a + 1)(a + 1) - 2 \sum_{k=0}^{a}k^{h}$ is the normalizing constant. It is symmetric around the target $x$, satisfying  Definition~\ref{Defkern} and suitable for count variables; see~\cite{KSKZ07}   and also \cite{Kokonendji&Zocchi10} for an asymmetric version. 
\item

From \cite{AJG76},~\cite{KSK11} deduced the following discrete kernel  that we here label DiracDU  (\texttt{DirDU}) as ``Dirac Discrete Uniform''. For  fixed $c \in \{2,3,\ldots\}$ the number of categories, we define $\mathbb{S}_{x,c} =\{0, 1, \ldots, c-1\}$ and  
\begin{equation*}
DU_{x,h;c}(u) =  \left( 1-h\right) \mathds{1}_{\left\{{x}\right\}}(u)+\dfrac{h}{c-1}\mathds{1}_{\mathbb{S}_{x,c}\setminus\left\{{x}\right\}}(u),
\end{equation*}
where $h \in (0,1]$ and $x \in \mathbb{T}_1$.  This DiracDU kernel is symmetric around the target, satisfying Definition~\ref{Defkern}  and appropriated for categorical set $\mathbb{T}_1$. See, e.g.,  \cite{Racine&Li04} for some uses.
%This kernel fullfills Definition~\ref{Defkern} and is adapted for categorical variables.  The DiracDU kernel was introduced in its multivariate version by  \cite{AJG76} and the current version by \cite{KSK11}; see \cite{Racine&Li04} for a use  in regression.
\item From the well known   \cite{E69} kernel $K^{E}(u) = \frac{3}{4}(1 - u^{2})\mathds{1}_{[ -1,1]}(u)$, we define its associated version  (\texttt{Epan}) on $\mathbb{S}_{x,h}= [ x-h,x+h]$ with $x \in\mathbb{T}_1:=\mathbb{R}$ and $h>0$:
\begin{equation*}
K^E_{x,h}(u)=\frac{3}{4h}\left\{1 - \left(\frac{u-x}{h}\right)^{2}\right\}\mathds{1}_{[ x-h,x+h]}(u). \label{gam2}
\end{equation*} 
It is obtained through (\ref{classcical}) and is well adapted for continuous variables with unbounded supports.
\item The gamma kernel (\texttt{Gamma}) is defined  on $\mathbb{S}_{x,h}= [0,\infty)=\mathbb{T}_1$ with  $x\in \mathbb{T}_1$ and  $h>0$: \begin{equation*}
GA_{x,h}(u)=\dfrac{u^{x/h}}{\Gamma\left(1+x/h\right)h^{1+x/h}}\exp{\left(-\dfrac{u}{h}\right)}
\mathds{1}_{[ 0,\infty)}(u), \label{gam2}
\end{equation*} 
where $\Gamma(\cdot)$ is the classical gamma function.  It is the pdf of the gamma distribution  $\mathcal{G}a(1 + x/h,h)$ with scale parameter $1 + x/h$ and shape parameter $h$.  It satisfies Definition~\ref{Defkern} and suits for non-negative  real set $\mathbb{T}_1$; see~\cite{Chen00a}.

\item  The  beta kernel (\texttt{Beta}) is however defined on $\mathbb{S}_{x,h}= [0,1]=\mathbb{T}_1$ with  $x\in \mathbb{T}_1$ and  $h>0$:
  \begin{equation*}
 BE_{x, h}(u) =\frac{u^{x/h}(1-u)^{(1 - x)/h}}{\mathscr{B}\left(1 + x/h, 1 + (1 - x)/h\right)} \mathds{1}_{[0, 1]}(u), \label{gam2}
\end{equation*} where $\mathscr{B}(r, s) = \int_{0}^{1}t^{r-1}(1 - t)^{s - 1}dt$ is the usual beta function with $r>0$ and $s>0$. It is the pdf of the beta distribution
 $\mathcal{B}e(1+x/h,(1-x)/h)$ with shape parameters $1+x/h$ and $(1-x)/h$. This pdf satisfies Definition~\ref{Defkern} and is appropriated for rates, proportions and percentages dataset $\mathbb{T}_1$; see~\cite{Chen99}.

\item We finally consider the bivariate beta kernel (\texttt{Bivariate beta}) defined by 
\begin{eqnarray}
\label{betakern}
 BS_{\mathbf{x}, \mathbf{H}}(u_1,u_2) &=& \left(\frac{u_{1}^{x_{1}/h_{11}}(1-u_{1})^{(1 - x_{1})/h_{11}} }{\mathscr{B}(1 + x_{1}/h_{11}, 1 + (1 - x_{1})/h_{11})}\right)  \left(\frac{u_{2}^{x_{2}/h_{22}}(1-u_{2})^{(1 - x_{2})/h_{22}}}{\mathscr{B}(1 + x_{2} / h_{22}, 1 + (1 - x_{2})/h_{22})} \right)  \nonumber \\
&& \times \left(1 + h_{12}\times\dfrac{u_{1} - \widetilde{\mu}_{1}(x_{1}, h_{11})}{h_{11}^{1/2}\widetilde{\sigma}_{1}(x_{1}, h_{11})}\times\dfrac{u_{2} - \widetilde{\mu}_{2}(x_{2}, h_{22})}{h_{22}^{1/2}\widetilde{\sigma}_{2}(x_{2}, h_{22})}\right)\mathds{1}_{\left[0, 1\right]^2}(u_1,u_2),
\end{eqnarray}
 with  $\mathbb{S}_{\mathbf{x}, \mathbf{H}}=\mathbb{T}_2=\left[0, 1\right]^2$, $\mathbf{x}=(x_1,x_2)^{\top} \in \mathbb{T}_2$ and $\mathbf{H} = \begin{pmatrix}h_{11}  & h_{12} \\ h_{12} & h_{22}\end{pmatrix}$.
For $j=1,2$,  the characteristics  in (\ref{betakern}) are given by  $h_{jj}>0$,  $\widetilde{\mu}_{j}(x_{j}, h_{jj}) = (x_{j} + h_{jj})/(1 + 2h_{jj})$,  $\widetilde{\sigma}_{j}^{2}(x_{j}, h_{jj}) = (x_{j} + h_{jj})(1 + h_{jj} - x_{j})(1 + 2h_{jj})^{-2}(1+3h_{jj})^{-1}h_{jj}$ , and the constraints 
\begin{equation}h_{12} \in \left[-\beta, \beta^{\prime}\right] \cap \left(-\sqrt{h_{11}h_{22}}\:,\sqrt{h_{11}h_{22}}\right) \label{constrainteSarmanov}
\end{equation} with 
$
\beta  = \left(\max_{v_1,v_2} \left\{\dfrac{v_{1} - \widetilde{\mu}_{1}(x_{1}, h_{11})}{h_{11}^{1/2}\widetilde{\sigma}_{1}(x_{1}, h_{11})}\times\dfrac{v_{2} - \widetilde{\mu}_{2}(x_{2}, h_{22})}{h_{22}^{1/2}\widetilde{\sigma}_{2}(x_{2}, h_{22})}\right\}\right)^{-1}
$ and \\
$
\beta^{\prime}  = \left\lvert\left(\min_{v_1,v_2} \left\{\dfrac{v_{1} - \widetilde{\mu}_{1}(x_{1}, h_{11})}{h_{11}^{1/2}\widetilde{\sigma}_{1}(x_{1}, h_{11})}\times\dfrac{v_{2} - \widetilde{\mu}_{2}(x_{2}, h_{22})}{h_{22}^{1/2}\widetilde{\sigma}_{2}(x_{2}, h_{22})}\right\}\right)^{-1}\right\lvert.
$
\newline
It satisfies Definition~\ref{Defkern} and is adapted for bivariate  rates.
The full bandwidth matrix $\mathbf{H}$ allows any  orientation of the kernel. Therefore, it can reach any point of the  space which might be inaccessible with  diagonal matrix. This type of kernel  is called beta-Sarmanov   kernel by \cite{Kokonendji&Some14}; see \cite{S66} and also \cite{L96} for this construction of multivariate densities with correlation structure from independent components. Like \cite{Bertin&Klutnitchkoff14}, the  miminax properties of this bivariate beta kernel are also possible and more generally for associated kernels.
\begin{figure}[!h]
\begin{center}
  \mbox{
\subfloat[(a) ]{\includegraphics[width=210pt,height=220pt,scale=0.95]{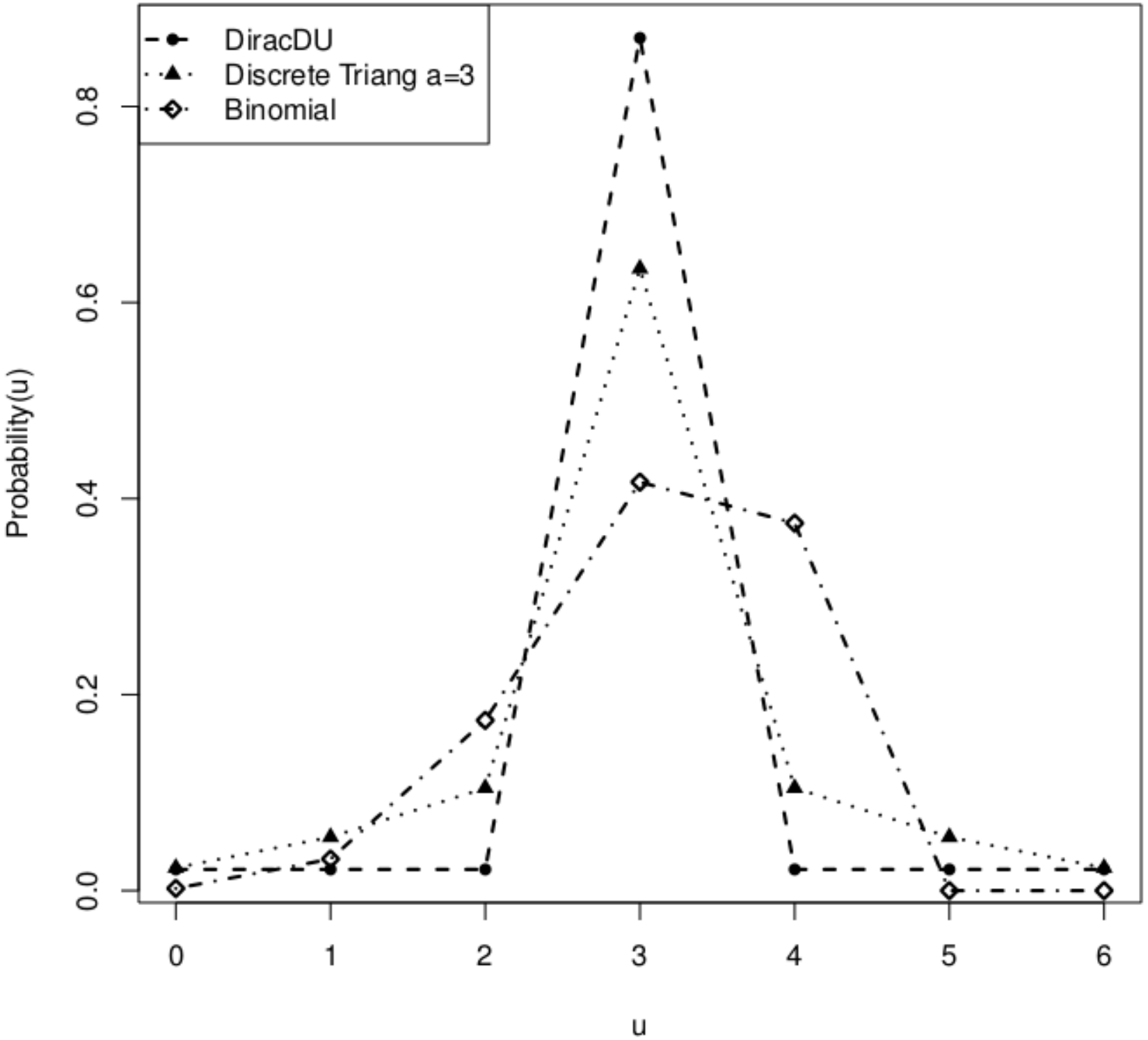} } \:
\subfloat[(b)]{\includegraphics[width=210pt,height=220pt,scale=0.95]{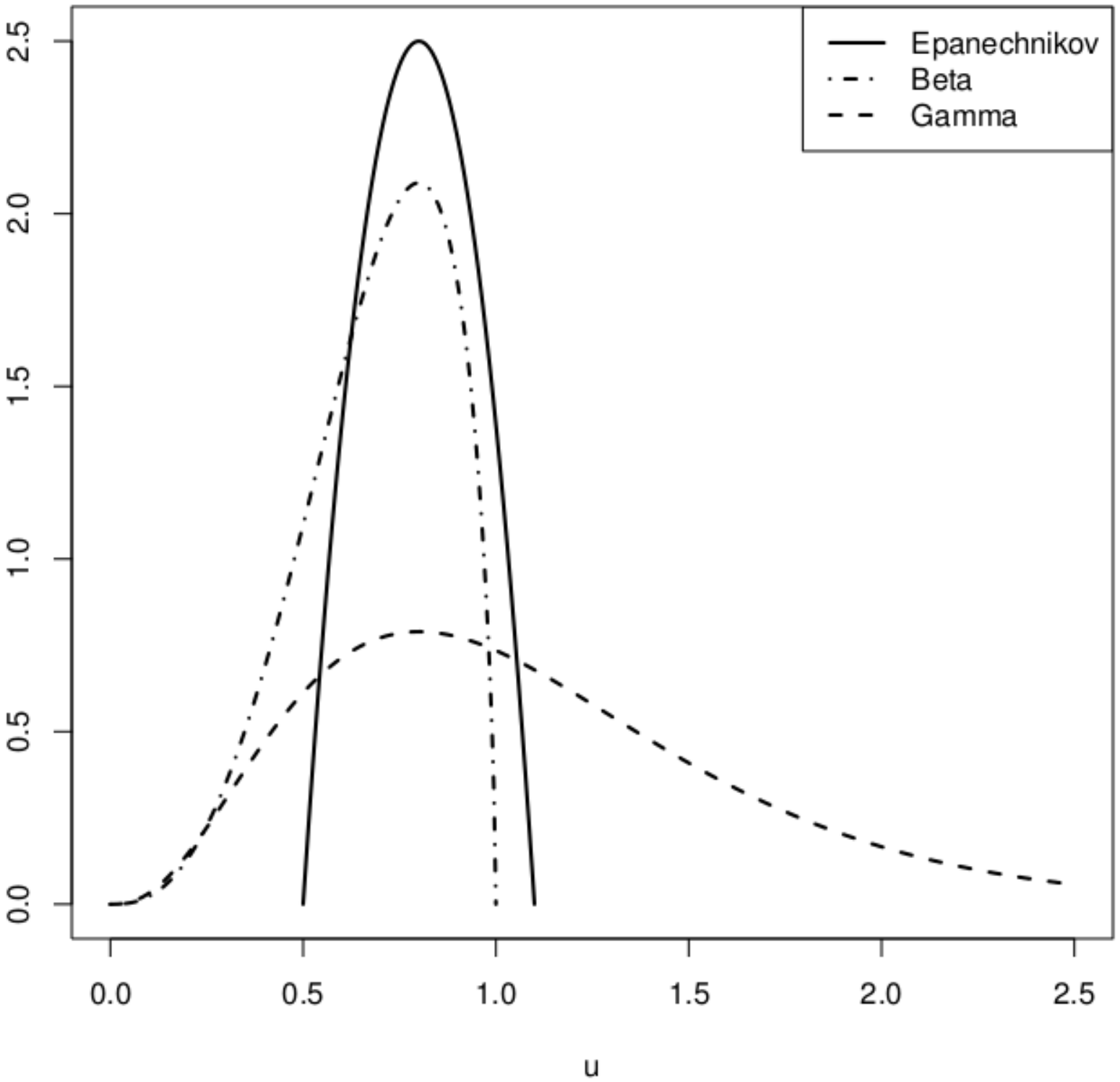} }  } 
 \end{center}
	\caption{Shapes of univariate (discrete and continuous) associated kernels: (a)  DiracDU,  discrete triangular $a=3$ and binomial with same target $x=4$ and bandwidth $h=0.13$; (b)   Epanechnikov, beta and gamma with same $x=0.8$ and $h=0.3$.}\label{Associatedkernels}
\end{figure}
\end{itemize}

Figure~\ref{Associatedkernels} shows some forms of the above-mentioned univariate associated kernels. The plots highlight the importance given to the target point and around it  in both discrete and continuous cases. Furthermore, for a fixed bandwidth 
 $h$,  the classical associated kernel of Epanechnikov, and also the categorical DiracDU kernel, keep their respective same shapes along the support; however, they  change according to the target for the others non-classical associated kernels. This explains the inappropriateness of the Epanechnikov kernel for density or regression estimation in any  bounded interval (Figure~\ref{Associatedkernels}(a)) and of the DiracDU kernel for count regression estimation (see simulations below).

\subsection{Bandwidth matrix selection by cross validation}
\label{ssec:bandwidth matrix selection}
In the context of multivariate kernel regression, the bandwidth matrix selection is here obtained by the well-known least squares cross-validation. In fact, for a given associated kernel, the optimal bandwidth matrix is
$\mathbf{\widehat{H}}  =  \displaystyle \mbox{arg }  \underset{\mathbf{H} \in \mathcal{H} }{ \min} \mbox{ LSCV}(\mathbf{H})$ with
\begin{equation}\label{Hcv}
  \mbox{LSCV}(\mathbf{H})  =  \frac{1}{n} \displaystyle\sum_{i = 1}^{n}\left\{Y_i - \widetilde{m}_{-i}(\mathbf{X}_{i}; \boldsymbol{\kappa})\right\}^{2},
\end{equation}
where $\widetilde{m}_{-i}(\mathbf{X}_{i}; \boldsymbol{\kappa})$ is computed as $\widetilde{m}_{n}$ of (\ref{RegFull}) excluding $\mathbf{X}_{i}$ and, $\mathcal{H}$ is the set of bandwidth matrices $\mathbf{H}$; see, e.g., \cite{Kokonendji&Senga&Demetrio09} in univariate case and also \cite{Zhangetal2014} and \cite{ZAK1a} for univariate bandwidth estimation by sampling algorithm methods. For diagonal bandwidth matrices (i.e. multiple associated kernels)  the LSCV method use the set of diagonal matrices $\mathcal{D}$. Concerning the beta-Sarmanov kernel (\ref{betakern}) with full bandwidth matrix, this LSCV method is used under $\mathcal{H}_1$, a subset of $\mathcal{H}$ verifying the constraint (\ref{constrainteSarmanov}) of the associated kernel. Their algorithms are described below and used for numerical studies in the following section.

\subsubsection*{\textit{Algorithms of LSCV method (\ref{Hcv}) for some type of associated kernels and their correponsding  bandwidth matrices}}\label{par:Algorithme}

\begin{enumerate}
\item[A1.] Bivariate beta (\ref{betakern}) with full bandwidth matrices and dimension $d=2$. 
%\begin{algorithm}
\begin{enumerate}[1.]
\item Choose two intervals $H_{11}$ and $H_{22}$ related to $h_{11}$ and $h_{22}$, respectively.
\item For $\delta = 1, \ldots, \ell(H_{11})$ and $\gamma = 1, \ldots, \ell(H_{22})$,
\begin{enumerate}[(a)]
 \item Compute the interval $H_{12}[\delta,\gamma]$ related to $h_{12}$ from constraints  in (\ref{constrainteSarmanov}); 
 \item For $\lambda = 1, \ldots, \ell(H_{12}[\delta,\gamma])$, \\
  Compose the full bandwidth matrix $\mathbf{H}(\delta,\gamma,\lambda):=\left(h_{ij}(\delta,\gamma,\lambda)\right)_{i,j=1,2}$ with $h_{11}(\delta,\gamma,\lambda)=H_{11}(\delta)$, $h_{22}(\delta,\gamma,\lambda)=H_{22}(\gamma)$ and 
  $h_{12}(\delta,\gamma,\lambda)=H_{12}[\delta,\gamma](\lambda)$.
   \end{enumerate}
  \item Apply LSCV method on the set $\mathcal{H}_1$ of all full bandwidth matrices $\mathbf{H}(\delta,\gamma,\lambda)$.
\end{enumerate}

%\item[A2.] Scott bandwidth matrices.  
%\begin{enumerate}[1.]
%\item Choose an interval $H$ related to $h$ and a fixed bandwidth matrix $\mathbf{H}_{0} = \left(h_{0ij}\right)_{i,j=1,2}$.
%\item For $\zeta = 1, \ldots, \ell(H)$, 
%\begin{enumerate}[(a)]
% \item Compute the interval $H_{012}[\zeta]$ related to $h_{012}$ from constraints (\ref{h12});
% \item For $\kappa = 1, \ldots, \ell(H_{012}[\zeta])$, \\
%  Compose the given bandwidth matrix $\mathbf{H}_{0}(\zeta,\kappa):=\left(h_{0ij}(\zeta,\kappa)\right)_{i,j=1,2}$ with $h_{011}(\zeta,\kappa)=h_{011}$, $h_{022}(\zeta,\kappa)=h_{022}$ and 
%  $h_{012}(\zeta,\kappa)=H_{012}[\zeta](\kappa)$;
%  \item Compose then the Scott bandwidth matrix $\mathbf{H}(\zeta,\kappa) := %H(\zeta)\times \mathbf{H}_{0}(\zeta,\kappa)$.
%   \end{enumerate}
%  \item Apply LSCV method on the set $\mathcal{S}_1$ of all Scott bandwidth matrices $\mathbf{H}(\zeta,\kappa)$.
%\end{enumerate}

\item[A2.] Multiple associated kernels (i.e. diagonal bandwidth matrices) for  $d \geq 2$. 
\begin{enumerate}[1.]
\item Choose two intervals $H_{11}$, $\ldots$, $H_{dd}$ related to $h_{11}$, $\ldots$, $h_{dd}$, respectively.
\item For $\delta_{1} = 1, \ldots, \ell(H_{11})$, $\ldots$, $\delta_{d} = 1, \ldots, \ell(H_{dd})$, \\
  Compose the diagonal bandwidth matrix $\mathbf{H}(\delta_{1},\ldots, \delta_{d}):=\mathbf{Diag}\left(H_{11}(\delta_1), \ldots, H_{dd}(\delta_d)\right)$.
\item Apply LSCV method on the set $\mathcal{D}$ of all diagonal bandwidth matrices $\mathbf{H}(\delta_1,\ldots, \delta_d)$.
\end{enumerate}
%\end{algorithm}
\end{enumerate}

 For a given interval $I$, the notation $\ell(I)$ is the total number of subdivisions of $I$ and $I(\eta)$ denotes the real value at the subdivision $\eta$ of $I$. Also, for practical uses of (A1) and (A2), the intervals $H_{11},\ldots, H_{dd}$ are taken generally according to the chosen associated kernel.

%We now examine the measure the effect of the type of associated kernel $\boldsymbol{\kappa}$ via 

\section{Simulation studies and real data analysis}
\label{sec:Simulation studies and real data analisys}
We apply the multivariate associated kernel estimators $\widetilde{m}_{n}$ of (\ref{RegFull}) and (\ref{Regprod}) to some simulated target regressions functions $m$  and then to two real datasets. The multivariate and multiple associated kernels used are built from those  of Section~\ref{ssec:Associated kernels for illustration}. The optimal bandwidth matrix is here chosen by LSCV method (\ref{Hcv}) using Algorithms A1 and  A2   of Section~\ref{par:Algorithme} and their indications. Besides the criterion of kernel support, we retain three measures to examine the effect of different associated kernels  $\boldsymbol{\kappa}$ on multiple regression. In simulations, it is the average squared errors (ASE)  defined as $$ASE(\boldsymbol{\kappa})=\dfrac{1}{n}\sum_{i=1}^{n}
\left\lbrace m(\mathbf{x}_i)-\widetilde{m}_{n}(\mathbf{x}_i;\boldsymbol{\kappa}) \right\rbrace^{2}.$$ For real datasets, we use the root mean squared error  (RMSE) which linked to ASE through squared root and by changing the simulated value $m(x_i)$ into the observed value $y_i$: 
\begin{equation*}
RMSE( \boldsymbol{\kappa}) = \sqrt{\frac{1}{n}\displaystyle\sum_{i=1}^n \left\{y_i -\widetilde{m}_{n}(\bold{x}_i;\boldsymbol{\kappa})\right\}^2 }.
\end{equation*} Also, we consider the practical  coefficient of determination $R^2$ which quantifies the proportion of variation of the response variable $Y_i$ explained by the non-intercept regressor $\bold{x}_i$
\begin{equation*} \label{R2}
R^2( \boldsymbol{\kappa}) = \frac{\sum_{i=1}^n \left\{\widetilde{m}_{n}(\bold{x}_i; \boldsymbol{\kappa}) - \overline{y}\right\}^2}{\sum_{i=1}^n (y_i - \overline{y})^2},
\end{equation*}
with $\overline{y} = n^{-1}(y_1 + \ldots + y_n)$. All these criteria above have their simulated or real data counterparts by replacing $y_i$ with $m(\bold{x_i})$ and vice versa. Computations have been performed on the supercomputer facilities of the M\'esocentre de calcul de Franche-Comt\'e using the R software; see \cite{R13}.
 \subsection{Simulation studies}
 \label{sec:Simulation studies}
 Expect as otherwise, each result is obtained with the number of replications $N_{sim}=100$.
  \subsubsection{Bivariate cases}
We consider seven   target regression functions  labelled A,  B, C, D and E with dimension $d=2$.

\begin{itemize}
 \item Function A is a bivariate beta without correlation $\rho(x_1,x_2) = 0$: 
\begin{equation*}\label{Betabiv}
    m(x_{1}, x_{2}) = \frac{x_1^{p_{1}-1}(1 - x_1)^{q_{1}-1}x_2^{p_{2}-1}(1 - x_2)^{q_{2}-1}}{\mathscr{B}(p_{1}, q_{1})\mathscr{B}(p_{2}, q_{2})} \mbox{  }  \mathds{1}_{\left[0, 1\right]}(x_{1})\mathds{1}_{\left[0, 1\right]}(x_{2}),
\end{equation*}
with $(p_{1}, q_{1}) = (3,2)$ and $(p_{2},q_{2}) = (5,2)$ as parameter values in univariate beta density.
\item  Function B is the bivariate Dirichlet density
\begin{equation*}\label{Dirichletbiv}
 m(x_{1}, x_{2}) = \frac{\Gamma(\alpha_{1} + \alpha_{2} + \alpha_{3})}{\Gamma(\alpha_{1})\Gamma(\alpha_{2})\Gamma( \alpha_{3})}x_{1}^{\alpha_{1} - 1}x_{2}^{\alpha_{2} - 1}(1 - x_{1} - x_{2})^{\alpha_{3} - 1} \mathds{1}_{\{x_{1},~ x_{2} \geq 0,~ x_{1} + x_{2} \leq 1\}}(x_{1},x_2),
\end{equation*}
where $\Gamma(\cdot)$ is the classical gamma function, with parameter values $\alpha_{1} = \alpha_{2} = 5$, $\alpha_{3} =6$ and, therefore, the moderate value of $\rho(x_1,x_2)=-(\alpha_{1}\alpha_{2})^{1/2}(\alpha_{1}+\alpha_{3})^{-1/2}(\alpha_{2}+\alpha_{3})^{-1/2}= -0.454$. 
%\item  density C is 
% \begin{equation}\label{CorContinuousbiv}
 %   m(\bold{x}) = sin(2\pi x_1) + x_2, \mathbb{   } \mbox{  }\mbox{  }  \bold{x} = (x_1, x_2)^{\top} \in \left[0, 1\right]\times\left[0, 1\right];
%\end{equation}
 \item Function C is a bivariate Poisson with null correlation $\rho(x_1,x_2) = 0$: 
\begin{equation*}\label{Poissonbiv}
    m(x_{1}, x_{2}) = \frac{e^{-5}2^{x_1}3^{x_2}}{x_1!x_2!} \mathds{1}_{\mathbb{N}}(x_{1})\mathds{1}_{ \mathbb{N}}(x_{2}).
\end{equation*}
\item  Function D is a bivariate Poisson with correlation structure
\begin{equation*}\label{PoissonbivCor}
 m(x_{1}, x_{2}) = e^{-(\theta_1+\theta_2+\theta_{12})} \displaystyle \sum_{i=0}^{min(x_1,x_2)}\frac{\theta_{1}^{x_1+i} \theta_{2}^{x_2+i} \theta_{12}^{i}}{(x_1+i)!(x_2+i)!i!}\mathds{1}_{\mathbb{N}\times \mathbb{N}}(x_{1}, x_{2}),
\end{equation*} with parameter  values $\theta_1 = 2$, $\theta_2 = 3$ and $\theta_{12} = 4$ and, therefore, the moderate value of $\rho(x_1,x_2)=\theta_{12}(\theta_{1}+\theta_{12})^{-1/2}(\theta_{2}+\theta_{12})^{-1/2}= 0.617$; see, e.g.,  \cite{YS}.
\item Function E is a bivariate beta without correlation $\rho(x_1,x_2) = 0$: 
\begin{equation*}\label{BetaPoissonbiv}
    m(x_{1}, x_{2}) = \frac{x_1^{p_{1}-1}(1 - x_1)^{q_{1}-1}3^{x_2}}{e^{3}\mathscr{B}(p_{1}, q_{1})x_2!} \mbox{  } \mathds{1}_{\left[0, 1\right]}(x_{1}) \mathds{1}_{ \mathbb{N}}(x_{2}),
\end{equation*}
with $(p_{1}, q_{1}) = (3,3)$.
%\item  density G 
\end{itemize}

\begin{table}[!h]
\begin{center}
\begin{tabular}{rrr}
\hline \hline
$n$&\multicolumn{1}{c}{Bivariate beta} & \multicolumn{1}{c}{Beta$\times$Beta} \\\hline
$50$ &$276.198$   &$7.551$  \\
$100$  &$647.255$& $30.081$ \\
\hline \hline
\end{tabular}
\caption{Typical Central Processing Unit (CPU) times (in seconds) for one replication  of LSCV method (\ref{Hcv}) by using Algorithms A1 and A2  of Section~\ref{ssec:bandwidth matrix selection}.} \label{Timehcv}
\end{center}
\end{table}
Table~\ref{Timehcv} presents the execution times needed for computing the LSCV method for both bivariate  beta kernels with respect to only one replication of sample sizes $n=50$  and $100$  for  the target function A. The computational times of the LSCV method for  the bivariate beta with correlation structure (\ref{betakern}) are obviously longer than those without correlation structure.  Let us note that for full bandwidth matrices, the execution times become very large when the number of observations is large; however, these CPU times  can be considerably reduced by parallelism processing, in particular for the bivariate beta kernel  with full LSCV method (\ref{Hcv}). These constraints (\ref{constrainteSarmanov}) reflect the difficulty for finding the appropriate bandwidth matrix with correlation structure by LSCV method. 

% where $N =100$ is the  number of replications. In particular $\boldsymbol{\nu}(d\mathbf{x}) = \nu_{1}(d{x}_1)\otimes \cdots  \otimes \nu_{d}(dx_d)$  where, for $j \in 1, \ldots, n$, $\nu_{j}$  is a Lesbesgue or count measure on the corresponding univariate support $\mathbb{T}^{[j]}_1$. In addition we use also $$\overline{R^2}= \dfrac{1}{N}\sum_{\ell=1}^{N}R^{2}_{\ell}(\boldsymbol{\kappa}),$$ where $R^{2}_{\ell}$ is defined in (\ref{R2}).
%The following results are given for only one simulation of (\ref{Betabiv}), (\ref{Dirichletbiv}), (\ref{Poissonbiv}) and  (\ref{PoissonbivCor}). %Figure~\ref{figReg} gives the regression function in black and the estimated one in grey.
\begin{table}[!h]
\begin{center}
\begin{tabular}{rrrrr}
\hline
\hline
&$n$&\multicolumn{1}{c}{Bivariate beta} & \multicolumn{1}{c}{Beta$\times$Beta}& \multicolumn{1}{c}{{\it Epan}$\times${\it Epan}} \\\hline

  \multirow{2}{*}{A}&   \multirow{1}{*}{50}%&   %$\begin{pmatrix}0.005& -1.940318\!\E\!-\!7 \\ -1.940318\!\E\!-\!7 & 0.005\end{pmatrix}$ & $\mathscr{D}\left(0.005,   0.006\right)$&$\mathscr{D}\left(0.315,   0.089\right)$ \\
   %& $0.9506(0.0332)$ &$0.9519(0.0314)$ &$0.5309(0.1213)$\\
       & $0.4368(0.3754)$&$0.4266(0.3724)$&$0.7483(0.2342)$\\
  &\multirow{1}{*}{100}%& $0.9565(0.0171)$ &$0.9579(0.0118)$  & \\ 
 &  $0.1727(0.0664)$&$0.1952(0.0816)$&0.6727(0.1413) \\
 \hline
   \multirow{2}{*}{B}& 50  & $1.2564(0.5875)$&$1.4267(0.4024)$&$1.6675(2.0353)$\\
 & 100   & $0.3041(0.1151)$&$0.3362(0.1042)$&$1.3975(1.5758)$\\
 % & $0.2186946$&$0.2186947$\\
\hline \hline
\end{tabular}\caption{Some expected values of 
$\overline{ASE}(\boldsymbol{\kappa})$ and their standard errors in  parentheses  with $N_{sim}=100$  of some  multiple associated kernel regressions for simulated continuous data  from functions A with $\rho(x_1,x_2) = 0$ and B with $\rho(x_1,x_2) = -0.454$.} \label{ErrBetabiv}
\end{center}
\end{table}
 
%\begin{table}[!h]
%\begin{center}
%\begin{tabular}{rcccccccc}
%\hline
%\hline
%$n$&\multicolumn{1}{c}{Bivariate beta} & \multicolumn{1}{c}{Beta$\times$Beta}& \multicolumn{1}{c}{{\it Epan}$\times${\it Epan}}  \\\hline

  % \multirow{1}{*}{50} % $\begin{pmatrix}0.007& -2.021594\!\E\!-\!7 \\ -2.021594\!\E\!-\!7 & 0.006\end{pmatrix}$ & $\mathscr{D}\left(0.007,0.004\right)$&$\mathscr{D}\left(0.231,   0.087\right)$ \\
   % & $0.9115(0.0738)$ &$0.9240(0.0318)$&$0.5459(0.1035)$\\
 %     & $1.2564(0.5875)$&$1.4267(0.4024)$&$1.6675(2.0353)$\\
 % \multirow{1}{*}{100} %&$\begin{pmatrix}0.003& -3.716716\!\E\!-\!8 \\ -3.716716\!\E\!-\!8 & 0.003\end{pmatrix}$  &$\mathscr{D}\left(0.003,0.003\right)$&$\mathscr{D}\left(0.259,   0.039\right)$  \\ 
 %&$0.9529(0.0184)$ &$0.9523(0.0194)$&$0.5113(0.1081)$\\
%  & $0.3041(0.1151)$&$0.3362(0.1042)$&$1.3975(1.5758)$\\
%\hline \hline
%\end{tabular}\caption{Some expected values of 
%$\overline{ASE}(\boldsymbol{\kappa})$ and their standard errors in  parentheses with $N_{sim}=100$  of some  multiple associated kernel regressions for simulated continuous data   from fucntion with $\rho(x_1,x_2) = -0.454$.} \label{ErrDirichletbiv}
%\end{center}
%\end{table} 
Table~\ref{ErrBetabiv}  reports the average $ASE(\boldsymbol{\kappa})$ which we denote $\overline{ASE}(\boldsymbol{\kappa})$ for
three continuous associated kernels $\boldsymbol{\kappa}$ with respect to functions A and B and according to sample sizes $n \in \{50, 100\}$.  We can see that  both  beta kernels in dimension $d=2$  work better than the multiple Epanechnikov  kernel for all sample sizes and all correlation structure in the regressors.  This reflects  the appropriateness of  the beta kernels which are suitable to the support of rate regressors.   Then, the explanatory variables with correlation structure  give larger $\overline{ASE}(\boldsymbol{\kappa})$  than those without correlation structure. Also, both  beta kernels  give quite similar results.
 Furthermore, all $\overline{ASE}(\boldsymbol{\kappa})$ are better when the sample size increases.

Finally, Tables~\ref{Timehcv} and \ref{ErrBetabiv}  highlight that the use of  bivariate beta kernels with correlation structure is not recommend in regression  with rates explanatory variables. Thus, we focus on multiple associated kernels for the rest of the simulations studies.
\begin{table}[!h]
\begin{center}
\begin{tabular}{rrccccccccc}
\hline
\hline
&$n$&\multicolumn{1}{c}{DTr2$\times$DTr2}&\multicolumn{1}{c}{DTr3$\times$DTr3 }&\multicolumn{1}{c}{Bin$\times$Bin}  & \multicolumn{1}{c}{{\it Epan}$\times${\it Epan}}  &\multicolumn{1} {c}{{\it DirDU}$\times${\it DirDU}} \\\hline

%  \multirow{3}{*}{20}   & (0.001, 0.005) &(0.001, 0.004)& (0.026, 0.321) & (1.001, 1.006)   &(0.005,0.029) \\
 %   & $[\mathbf{98.476}]$ &$[97.928]$&$[59.101]$&$[ 98.065]$&$[99.971]$\\
   \multirow{3}{*}{C}& \multirow{1}{*}{20}   &1.5e-6(2.2e-6)&3.3e-6(4.1e-6)&3.6e-5(9.7e-6)&4.0e-5(3.5e-5)&1.6e-8(1.8e-8)\\ 
&\multirow{1}{*}{50} & 3.1e-7(6.9e-7) &4.7e-7(9.7e-7)& 3.6e-5(7.4e-6) &3.8e-5(2.8e-5) &3.7e-9(2.3e-9) \\
    
   %   &$[\mathbf{99.545}]$ &$[99.167]$&$[54.346]$& $[\mathbf{99.566}]$&$[99.996]$\\
    %    &$7.957$ &$\underline{1.854}$ & $30.221$&$2.136$  \\\\
   
   & \multirow{1}{*}{100}&8.6e-8(1.2e-7)&2.9e-7(3.1e-7)& 3.7e-5(4.8e-6)&3.6e-5(2.3e-5)& 4.1e-10(3.5e-10) \\
   %&$[\mathbf{99.548}]$ &$[99.212]$ & $[47.383]$& $[66.531]$&$[100]$\\
  % &$\underline{5.089}$ &$11.774$ & $10.273$&$6.470$ \\
\hline 
  \multirow{3}{*}{D}& \multirow{1}{*}{20}  & 2.4e-6(2.8e-6)&4.5e-6(4.9e-6) & 7.1e-6(2.6e-6) & 4.2e-6(2.5e-6)   &2.7e-8(2.1e-8) \\
  % &  $[69.515]$&$[\mathbf{75.800}]$& $[67.229]$& $[69.081]$&$[93.538]$\\
   % &$\underline{5.612}$&$12.096$&$27.426$&$3.988$\\\\
   
      & \multirow{1}{*}{50}&2.5e-7(3.4e-7)&1.8e-7(2.5e-7) & 8.1e-5(4.3e-6) & 5.1e-6(1.2e-6)   &4.3e-9(3.2e-9)\\
     % &  $[\mathbf{99.829}]$&$[99.726]$&$[71.840]$&$[82.734]$&$[99.993]$\\
     % &$\underline{4.211}$&$9.000$&$11.636$&$3.686$
          
    &\multirow{1}{*}{100}&2.6e-8(6.2e-8)&4.8e-8(9.5e-8)&9.3e-6(8.2e-7)&7.2e-6(7.8e-7)&5.3e-10(4.6e-10)\\
   %  &$[\mathbf{99.776}]$& $[98.585]$&&$[76.429]$&$[99.999]$\\
   %  &$\underline{1.890}$&$4.348 $&&$1.897$& \\
\hline \hline
\end{tabular}
\caption{ Some expected values of 
$\overline{ASE}(\boldsymbol{\kappa})$ and their standard errors in  parentheses with $N_{sim}=100$ of some  multiple associated kernel regressions for simulated count data  from functions C with $\rho(x_1,x_2) = 0$ and D with $\rho(x_1,x_2) = 0.617$.} \label{ErrPoissonbiv}
\end{center}
\end{table} 

%\begin{table}[!h]
%\begin{center}
%\begin{tabular}{rcccccccc}
%\hline
%\hline
%$n$&\multicolumn{1}{c}{DTr2$\times$DTr2}&\multicolumn{1}{c}{DTr3$\times$DTr3 }&\multicolumn{1}{c}{Bin$\times$Bin}  & \multicolumn{1}{c}{{\it Epan}$\times${\it Epan}}  &\multicolumn{1}{c}{{\it DirDU}$\times${\it DirDU}} \\\hline

%  \multirow{1}{*}{20}  & 2.4e-6(2.8e-6)&4.5e-6(4.9e-6) & 7.1e-6(2.6e-6) & 4.2e-6(2.5e-6)   &2.7e-8(2.1e-8) \\
  % &  $[69.515]$&$[\mathbf{75.800}]$& $[67.229]$& $[69.081]$&$[93.538]$\\
   % &$\underline{5.612}$&$12.096$&$27.426$&$3.988$\\\\
   
   %    \multirow{1}{*}{50}&2.5e-7(3.4e-7)&1.8e-7(2.5e-7) & 8.1e-5(4.3e-6) & 5.1e-6(1.2e-6)   &4.3e-9(3.2e-9)\\
     % &  $[\mathbf{99.829}]$&$[99.726]$&$[71.840]$&$[82.734]$&$[99.993]$\\
     % &$\underline{4.211}$&$9.000$&$11.636$&$3.686$
          
  %  \multirow{1}{*}{100}&2.6e-8(6.2e-8)&4.8e-8(9.5e-8)&9.3e-6(8.2e-7)&7.2e-6(7.8e-7)&5.3e-10(4.6e-10)\\
   %  &$[\mathbf{99.776}]$& $[98.585]$&&$[76.429]$&$[99.999]$\\
   %  &$\underline{1.890}$&$4.348 $&&$1.897$& \\
%\hline \hline
%\end{tabular}
%\caption{Some expected values of 
%$\overline{ASE}(\boldsymbol{\kappa})$ and their standard errors in  parentheses with $N_{sim}=100$ of some  multiple associated kernel regressions of simulated count data  from function D with $\rho(x_1,x_2) = 0.617$.} \label{ErrPoissonbivCor}
%\end{center}
%\end{table} 

Table~\ref{ErrPoissonbiv}  shows the values $\overline{ASE}(\boldsymbol{\kappa})$  with respect to five associated kernels $\boldsymbol{\kappa}$  for sample size $n=20, 50$ and $100$ and count datasets generated from C and D. Globally, the discrete associated kernels in multiple case  perform better than the multiple Epanechnikov kernel  for all sample sizes and correlation structure in the regressors.  The use of categorical DiracDU kernels gives the best result in term of $\overline{ASE}(\boldsymbol{\kappa})$ but  DiracDU does not suit for these count datasets. Also, the discrete triangular kernels gives the most interesting result with an advantage to the discrete triangular with small arm $a=2$. This discrete triangular is the best since it concentrates always on the target and a few observations around it; see  Figure~\ref{Associatedkernels}(a). The results become much better when the sample size increases. The values $\overline{ASE}(\boldsymbol{\kappa})$ for regressors with or without correlation structure are comparable; and thus, we can focus on target regression functions without correlation structure for the remaining simulations.

%The values $\overline{ASE}(\boldsymbol{\kappa})$ without correlation structure between the  regressors are better  than those with strong correlation. Finally, 

\begin{table}[!h]
\begin{center}
\begin{tabular}{rrcccccccc}
\hline
\hline
&$n$  &\multicolumn{1}{c}{Beta$\times$DTr2}&\multicolumn{1}{c}{Beta$\times$	DTr3}&\multicolumn{1}{c}{Beta$\times$Bin}  & \multicolumn{1}{c}{{Beta$\times$\it Epan}} &\multicolumn{1}{c}{Beta$\times${\it DirDU}} \\\hline

   \multirow{3}{*}{E}&\multirow{1}{*}{30}  &3.738(1.883) &1.966(1.382) & 3.884(1.298) & 6.361(2.134)&0.162(0.201)  \\ 
    % &  $[64.346]$&$[\mathbf{80.721}]$& $[43.258]$& $[28.412]$\\\\ \\
%    &$5.612$&$12.096$&$27.426$&$3.988$&$41.515$\\\\
   
       &\multirow{1}{*}{50}&3.978(1.404)&2.106(1.119) &3.683(0.833)&7.143(1.732)&0.138(0.171) \\
     % &  $[54.760]$&$[\mathbf{81.528}]$&$[31.178]$&
      
           &  \multirow{1}{*}{100}&3.951(1.052)&1.956(0.806) & 3.835(0.834)&7.277(1.574) &0.113(0.147)\\
\hline \hline
\end{tabular}
\caption{Some expected values  ($\times 10^3$) of 
$\overline{ASE}(\boldsymbol{\kappa})$ and their standard errors in  parentheses with $N_{sim}=100$ of some  multiple associated kernel regressions of simulated mixed data from function E with $\rho(x_1,x_2) = 0$. }\label{ErrBetaPoissonbiv}
\end{center}
\end{table}

Table~\ref{ErrBetaPoissonbiv}  presents the values  for sample sizes $n \in \{30, 50, 100\}$ and for five  associated kernels $\boldsymbol{\kappa}$. The datasets are  generated from E and the beta kernel is applied on the continuous rate variable of E. We observe the superiority of the  multiple associated kernels using discrete kernels  over those defined with the Epanechnikov kernel for all sample sizes. Then, the multiple associated kernel with the categorical DiracDU gives the best $\overline{ASE}(\boldsymbol{\kappa})$ but it  is not appropriate for the count variable of E. Also, the values $\overline{ASE}(\boldsymbol{\kappa})$ are getting better when the sample size increases. 

From Tables~\ref{ErrBetabiv}, \ref{ErrPoissonbiv} and \ref{ErrBetaPoissonbiv},  the importance  of the type of associated kernel $\boldsymbol{\kappa}$ which respect the support of the explanatory variables is proven. %Thus, we focus in the multivariate case only on appropriated type of kernels $\boldsymbol{\kappa}$.

  \subsubsection{Multivariate cases}
  Since the appropriate associated kernels perform better than the inappropriate ones, we focus in higher dimension $d>2$ on regression with only suitable associated kernels.
Then, we consider two target regression functions labelled F and G for $d=3$ and 4 respectively. The formulas of the functions are given below.
\begin{itemize}
 \item Function F is a 3-variate with null correlation: 
\begin{equation*}\label{BetaandPoissonbiv}
    m(x_{1}, x_{2},x_{3}) = \frac{x_1^{p_{1}-1}(1 - x_1)^{q_{1}-1}2^{x_2}3^{x_3}}{e^{5}\mathscr{B}(p_{1}, q_{1})x_2!x_3!}\mathds{1}_{\left[0, 1\right]}(x_{1})\mathds{1}_{\mathbb{N}}(x_{2})\mathds{1}_{\mathbb{N}}(x_{3}),
\end{equation*}
with $(p_{1}, q_{1}) = (3,2)$.
 \item Function G is a 4-variate without correlation: 
\begin{equation*}\label{BetabivandPoissonbiv}
    m(x_{1}, x_{2},x_{3},x_{4}) = \frac{x_1^{p_{1}-1}(1 - x_1)^{q_{1}-1}x_2^{p_{2}-1}(1 - x_2)^{q_{2}-1}2^{x_3}3^{x_4}}{e^{5}\mathscr{B}(p_{1}, q_{1})\mathscr{B}(p_{2}, q_{2})x_3!x_4!}\mathds{1}_{\left[0, 1\right]}(x_{1})\mathds{1}_{\left[0, 1\right]}(x_{2}) \mathds{1}_{\mathbb{N}}(x_{3})\mathds{1}_{\mathbb{N}}(x_{4}),
\end{equation*}
with $(p_{1}, q_{1}) = (3,2)$ and $(p_{2},q_{2}) = (5,2)$.
\end{itemize}

\begin{table}[!h]
\begin{center}
\begin{tabular}{rccccccccc}
\hline
\hline
$n$&\multicolumn{1}{c}{Beta$\times$DTr3$\times$DTr3}&\multicolumn{1}{c}{Beta$\times$Bin$\times$DTr3}&\multicolumn{1}{c}{Beta$\times$Beta$\times$DTr3$\times$DTr3} \\\hline
  \multirow{1}{*}{30}   &0.2501(0.1264)  & 0.3038(0.1258)&0.7448(0.5481)\\
  
   % & $[99.897]$ &\\
   % &\\ \\
    \multirow{1}{*}{50}&0.2381(0.0661) &0.2895(0.0162)&0.6055(0.2291)\\

     %  & & \\\\
   %
    \multirow{1}{*}{100}&0.2282(0.0649)&0.2822(0.0608)&0.5012(0.2166)\\
  % & & \\
\hline \hline%
\end{tabular}
\caption{Some expected values  ($\times 10^3$) of 
$\overline{ASE}(\boldsymbol{\kappa})$ and their standard errors in  parentheses with $N_{sim}=100$ of some  multiple associated kernel regressions of simulated mixed data  from 3-variate F and 4-variate G.} \label{ErrBetaandPoissonbiv}
\end{center}
\end{table} 
Table~\ref{ErrBetaandPoissonbiv} presents the regression study for dimension $d=3$ and 4  with respect to functions F and G and for sample size $n \in \{30, 50, 100\}$. The values $\overline{ASE}(\boldsymbol{\kappa})$  show the superiority of the multiple associated kernels using the discrete triangular kernel with $a=3$ over the one  with the binomial kernel. Some results with respect to function G for  an associated kernel  $\boldsymbol{\kappa}$ composed by two beta and two discrete triangular kernels with $a=3$ are also provided. The errors become smaller when the sample size increases.

\subsection{Real data analysis}
The dataset consists on a sample of 38 family  economies from a US large city and is available as the {\it FoodExpenditure} object in the {\it betareg} package of \cite{Cribari&Neta(2010)}. The dataset in its current form gives not available (NA) responses for associated kernel regressions especially when we use the discrete triangular or the  DiracDU kernel. Then, we extend  the original {\it FoodExpenditure} dataset with its first 20  observations which guarantees some results for the regression, and thus $n=58$. The dependent variable is {\it food/income}, the proportion of household \textit{income} spent on \textit{food}. Two explanatory variables are available:  the previously mentioned household {\it income} ($x_1$) and the {\it number of residents} ($x_2$) living in the household with $\widehat{\rho}(x_1,x_2) = 0.028$. We use the Gamma or the Epanechnikov kernel for the  continuous variable $income$  and the discrete (of Figure~\ref{Associatedkernels}(a))  or the Epanechnikov for the count variable {\it number of residents}. 

The results of the multiple associated kernels for regression are divided in two in Table~\ref{dataset1}. The appropriate  associated kernels which strictly follow the support of each variable give comparable results in terms of both RMSE$(\boldsymbol{\kappa})$ and R$^2(\boldsymbol{\kappa})$. In fact, the associated kernels that use the discrete triangular with arm $a=2$ and 3 give some R$^2(\boldsymbol{\kappa})$ approximately equal to $64 \%$. The inappropriate kernels give various results. The multiple Epanechnikov kernel and the type of kernel with DiracDU give  R$^2(\boldsymbol{\kappa})$ higher than $80 \%$ while the Gamma$\times$Epanechnikov gives R$^2(\boldsymbol{\kappa})$ less than $50\%$. Then, a little difference in terms of RMSE$(\boldsymbol{\kappa})$ can induce a high incidence on the R$^2(\boldsymbol{\kappa})$.

%Using kernels which follows scrupulously the support of the explanatory variables, Table~\ref{dataset1} shows the results of the nonparametric regressions. The   Gamma$\times$DTr2 kernel  with $a=2$ and  
%represents the most interesting results among these appropriate kernel regressions. Then, we have the Gamma$\times$DTr3  with $a=3$ and the Gamma$\times$Bin . 
%(RMSE$(\boldsymbol{\kappa})$=0.01426) (RMSE$(\boldsymbol{\kappa})$=0.01730)
% RMSE$(\boldsymbol{\kappa})$=0.01409 (RMSE$(\boldsymbol{\kappa})$=0.01278)
%(RMSE$(\boldsymbol{\kappa})$=0.01451) (RMSE$(\boldsymbol{\kappa})$=0.03266)

%In addition Table~\ref{dataset1} gives the same results with kernels whose supports  violate those of the explanatory variables. The  associated kernels  used gives varying results. The  Gamma$\times${\it DirDU}  and the   {\it Epan}$\times${\it Epan}  gives comparable results to the previous associated kernels while  Gamma$\times{\it Epan}$  gives much less. Finally, the choice of appropriate associated kernel according to the support of the explanatory variables is much more important than the values of RMSE$(\boldsymbol{\kappa})$ and R$^2(\boldsymbol{\kappa})$.

\begin{table}[!ht]
\begin{center}
\begin{tabular}{cccccccc}
\hline
\hline
  \multirow{3}{*}{ {Appropriate }}  &\multicolumn{1}{c}{Gamma$\times$DTr2}&\multicolumn{1}{c}{Gamma$\times$DTr3}&\multicolumn{1}{c}{Gamma$\times$Bin}  \\%\hline

%&$(0.271,0.001)$ &$(0.286,0.001)$ & $(0.061,0.211)$   \\ 
    &0.01409  &0.01426 &0.01730 \\
  &64.2681  &64.2708  &56.1091 \\ 
\hline %\hline
 \multirow{3}{*}{ {{\it Inappropriate }}} 	&\multirow{1}{*}{ {\it Epan}$\times${\it Epan}}&\multirow{1}{*}{ Gamma$\times${\it Epan}} & \multirow{1}{*}{ Gamma$\times${\it DirDU} } &
	 \\%\hline

 % & $(0.196,0.001)$ &$(1.334,2.061)$&$(0.211,0.001)$   \\ 
 &0.01451&0.03266 &0.01278\\
    &86.0011&47.0181 &89.3462\\

\hline \hline
\end{tabular}
\caption{Some expected values of RMSE$(\boldsymbol{\kappa})$ and in percentages R$^2(\boldsymbol{\kappa})$  of some multiple associated kernel regressions for the $FoodExpenditure$ dataset  with $n=58$.} \label{dataset1}
\end{center}
\end{table}
\begin{table}[!ht]
\begin{center}
\begin{tabular}{rrr|rrr|rrr}
\hline
\hline
   \multicolumn{1}{c}{$x_{1i}$}& \multicolumn{1}{c}{ $x_{2i}$}&\multicolumn{1}{c|}{$y_i$}& \multicolumn{1}{c}{$x_{1i}$}&  \multicolumn{1}{c}{$x_{2i}$}&\multicolumn{1}{c|}{$y_i$}& \multicolumn{1}{c}{$x_{1i}$}& \multicolumn{1}{c}{ $x_{2i}$}& \multicolumn{1}{c}{$y_i$}\\
  \hline
  68.1& 54.6&  0.8 & 80.3 &9.1 & 1.5  &78.6& 31.0 & 1.6\\
  60.8& 4.4&  1.3 &23.1 &83.1&  2.3 &44.9& 3.2 & 1.0\\
  34.4& 36.2&  1.2 &16.9 &90.4 & 2.0 &78.2& 13.9 & 1.8\\
  59.4& 27.5&  1.3 &9.4 &79.2 & 2.8 &60.2& 35.2 & 1.1\\
  4.7& 81.0&  2.9 &55.8 &21.9 & 1.3 &65.6& 26.1 & 1.6\\
  19.9& 97.4&  1.2 &27.5 &75.0 & 2.2  &74.4& 12.6& 1.6\\
  20.6& 73.6&  2.4 &59.1&12.9 & 1.4 &83.5& 13.3 & 1.8\\
  16.4 &42.9&  1.1 &2.7 &93.9 & 2.4 &10.9& 83.5 & 2.6\\
  29.9 &74.4&  2.0 &13.9& 56.9& 1.4 &27.0& 77.1 & 2.2\\
 84.8& 26.6&  1.6 &14.0& 92.9& 2.1&3.1& 67.0 & 2.2\\
 46.1 &66.9&  1.2 &22.9&43.9  &1.1  &14.8& 72.9 & 2.5\\
 10.2& 86.3&  2.5 &53.8 &56.2  &1.0 &80.6& 16.5&  1.6\\
 89.4 &32.5&  1.6 &23.7& 61.5 &1.5 &64.1& 28.6 & 1.5\\
 30.9& 46.3&  1.1&39.6 &67.2 & 1.4 &15.6& 90.5 & 2.0\\
 24.3 &37.8&  1.2 &59.5& 45.1 & 0.9 &3.9& 68.6 & 2.5\\
 27.4 &74.6 & 1.9 &17.3& 81.2 & 2.6 &66.9& 43.7& 0.9\\
 47.7 &61.7 & 1.1 &93.7& 28.5 & 1.5 &1.5& 65.8 & 2.3\\
 33.1 &83.8&  1.5 &28.7& 82.7 & 2.0 &35.6& 43.7 & 1.0\\
 0.3 &83.3&  3.0 &61.3& 70.9 & 0.6  &13.9& 25.0 & 0.8\\
 76.9& 35.4&  1.2 &67.1& 24.0 & 1.7 &13.2& 70.8 & 2.2\\
 29.5& 44.6&  1.3 &85.8& 36.5 & 1.2 &34.5& 73.7 & 1.8\\
 19.6 &67.7&  1.9 &35.5& 76.9 & 1.8 &55.6& 6.9 & 1.3\\
 96.2 &26.1&  1.7 &18.8& 55.9 & 1.3 &30.7& 9.1 & 0.9\\
 85.9& 28.0 & 1.5 &50.4& 17.7 & 1.4 &43.5&15.1 & 1.0\\
 5.6& 39.1 & 1.1 &67.2& 8.7 & 1.5&31.5& 36.7 & 1.2\\
 99.9 &7.15 & 1.3 &13.1& 59.4 & 1.7&30.0& 21.5 & 0.8\\
 61.0 &31.1 & 1.4 &13.7& 75.8 & 2.5\\
% 80.3 &9.1 & 1.5 & 44.9& 3.2 & 1.0\\
\hline \hline
\end{tabular}\caption{Proportions (in $\%$) of folks who like the company,  those who like its strong product and turnover  of a company,  designed respectively by the variables $x_{1i}$, $x_{2i}$ and $y_i$,  with $\widehat{\rho}(x_1,x_2) = -0.6949$ and $n=80$.} \label{data2}
\end{center}
\end{table}  
Table~\ref{data2} of the second dataset  aims to explain the turnover of a large company by two proportions explanatory variables obtained by survey. The first variable $x_1$ is the rate of people who like the company and the second one $x_2$ is the percentage of people who like the strong product of this company. The dataset is obtained in 80 branch of this company. Obviously, there is a significant correlation between these explanatory variables: $\widehat{\rho}(x_1,x_2)=-0.6949$.

Table~\ref{dataset2} presents the  results for the nonparametric regressions with three associated kernels $\boldsymbol{\kappa}$. Both  beta kernels offer the most interesting results with R$^2(\boldsymbol{\kappa})$ approximately equal to $86\%$. Note that, the multiple Epanechnikov kernel gives lower performance  mainly because this continuous unbounded kernel does not suit for these bounded explanatory variables.

% ($RMSE(\boldsymbol{\kappa}) =0.10524$)($RMSE(\boldsymbol{\kappa}) =0.18886$)In order to smooth this dataset, we apply both general and multiple beta kernel in Section~\ref{sec:Simulation studies}($RMSE(\boldsymbol{\kappa}) =0.10523$)

\begin{table}[!ht]
\begin{center}
\begin{tabular}{rcccccccc}
\hline
\hline
&\multicolumn{1}{c}{Bivariate beta} & \multicolumn{1}{c}{Beta$\times$Beta}& \multicolumn{1}{c}{ {\it Epan}$\times${\it Epan}} \\\hline

  \multirow{2}{*}{}%&   $\begin{pmatrix}0.019& 4.38847\!\E\!-\!6 \\ 4.38847\!\E\!-\!6 & 0.014\end{pmatrix}$ & $\mathscr{D}\left(0.019,   0.0.014\right)$&$\mathscr{D}\left(0.151,   0.223\right)$ \\
       & $0.10524$ & $0.10523$&$0.18886$\\
    & $86.6875$ & $86.6874$&$76.3431$\\
      
\hline \hline
\end{tabular}\caption{Some expected values of RMSE$(\boldsymbol{\kappa})$ and  in percentages R$^2(\boldsymbol{\kappa})$  of some bivariate associated kernel regressions for tunover dataset in Table~\ref{data2} with $\widehat{\rho}(x_1,x_2) = -0.6949$ and $n=80$.} \label{dataset2}
\end{center}
\end{table}

\section{Summary and final remarks}
\label{sec:Summary and final remarks}
%Certaines propri\'et\'es asymptotiques de ces estimateurs (\ref{RegFull}) et (\ref{Regprod}) seront pr\'esent\'ees en utilisant les d\'eriv\'ees et les diff\'erences finies. \`A travers l'erreur moyenne quadratique  et le coefficient de d\'etermination, des \'etudes num\'eriques (par simulation et sur donn\'ees r\'eelles) dans le cas bivari\'e montrent l'efficacit\'e de cette m\'ethode.   Le choix de la matrice des fen\^etres optimale se fera Ã  chaque fois par validation crois\'ee avant une am\'elioration par des approches bay\'esiennes. Ce choix implique naturellement des calculs fastidieux~;  notamment, dans le cas g\'en\'eral (\ref{RegFull}) avec structure de corr\'elation oÃ¹ l'on a Ã  s\'electionner une matrice sym\'etrique Ã  $d(d+1)/2$ paramÃ¨tres.

We have presented associated kernels for nonparametric multiple  regression and in presence of a mixture of discrete and continuous
explanatory variables;  see, e.g.,  \cite{ZAK14b} for a  choice of the bandwidth matrix by Bayesian methods.  Two particular cases including the continuous classical and the multiple (or product of) associated kernels are highlight  with the bandwidth matrix selection by cross-validation. Also, six univariate  associated kernels and a bivariate beta with correlation structure are presented and used for computational studies. 

Simulation experiments and analysis of two real datasets provide insight into the behaviour of the type of associated kernel $\boldsymbol{\kappa}$ for small and moderate sample sizes. Tables \ref{Timehcv}, \ref{ErrBetabiv} and  \ref{dataset2} on bivariate rate regressions can be conceptually summarized as follows. The use of associated kernels with correlation structure is not recommend. In fact, it is time consuming and have the same performance as the multiple beta kernel. Also, these appropriate beta kernels are better than the inappropriate multiple Epanechnikov. For count regressions, the multiple associated kernels built from the binomial and the discrete triangular  with small arms are superior to  those with the optimal continuous Epanechnikov. Furthermore, the categorical DiracDU kernel gives misleading results since it does not suit for count variables, see Tables~\ref{ErrPoissonbiv} and  \ref{ErrBetaPoissonbiv}.  We advise beta kernels for rates variables and gamma kernels for non-negative  dataset for small and moderate sample sizes, and also for all dimension $d \geq 2$; see, e.g., Tables~\ref{ErrBetaandPoissonbiv} and \ref{dataset1}. Finally, more than the performance of the regression, it is the correct choice of  the associated kernel according to the explanatory variables which is the most important. In other words, the criterion for choosing an associated kernel is  the support; however, for several kernels matching the support, we use common measures such as the mean integrated squared error. It should be noted that a large coefficient of determination $R^2$ does not mean good adjustment of the data; see Tables~\ref{dataset1} and \ref{dataset2}. Further research on associated kernels for functional regression is conceivable; see, e.g.,  \cite{ACT14} for  classical  kernels.

%Conclusion (partial):
%\begin{itemize}

 %\item No correlation effect for regressors in \textit{continuous} associated kernels (general/multiple)

 %\item No correlation effect for regressors in \textit{discrete} associated kernels (multiple)

 %\item Some correlation effects for regressors in \textit{mixed} associated kernels (multiple)
 %\item the results of simulations and real data analysis were predictable based on the shapes of the associated kernels in Figure~\ref{Associatedkernels}.

%\item<1-> R\'eduction du biais (9 r\'egions): \textit{Ã  illustrer}
%\end{itemize}

%Forthcoming works:
%\begin{itemize}
 % \item Improve bandwidth matrices choices (e.g. Bayesian approach)

%\item Construction of discrete/mixed associated kernels with correlation

%\item Simulation studies  with   regressors $>2$

%\item Application: Associated kernel discriminant analysis
%\end{itemize}

%\section*{Acknowledgements}
%We  sincerely thank Tristan Senga Ki\'ess\'e for useful discussions.

\section*{References}
%\bibliographystyle{elsarticle-num}
%\bibliographystyle{model1-num-names}
%\addcontentsline{toc}{chapter}{References}

%Further research of this work are in progress, especially on the discriminant analysis using the multiple associated kernels which have diagonal bandwidth matrices. Also, it will be interesting to investigate the choice of the bandwidth matrix by Bayesian methods; see, e.g.,  \cite{ZAK14b}. Furthermore, a construction and an implementation of appropriated associated kernels with correlation structure for discrete and mixed regressors  can extend the current work; see, e.g., \cite{S66} and \cite{L96}. Finally, we can use associated kernels for functional regression; see, e.g., \cite{ACT14} for  classical  kernels.%There are no correlation effect of the regressors in continuous associated kernels with or without correlation structure for regression. 

%\newpage
%\tableofcontents
\end{document}